\documentclass[11pt]{article}
\usepackage{amsmath}
\usepackage{amssymb}
\usepackage{amsfonts}
\usepackage[all,2cell]{xy} \UseAllTwocells \SilentMatrices
\usepackage{geometry,latexsym,amsmath,amstext,mathrsfs,amssymb,amsthm,amscd,amsopn,verbatim,graphics}
\usepackage[english]{babel}

\newtheorem{theorem}{Theorem}[section]
\newtheorem{proposition}[theorem]{Proposition}

\newtheorem{corollary}[theorem]{Corollary}
\newtheorem{definition}[theorem]{Definition}
\newtheorem{example}[theorem]{Example}
\newtheorem{examples}[theorem]{Examples}
\newtheorem{remark}[theorem]{Remark}
\numberwithin{equation}{subsection}

\newcommand*\Hom{\mathrm{Hom}}

\def\BA{{\mathbb{A}}}
\def\BC{{\mathbb{C}}}

\def\BK{{\mathbb{K}}}

\def\U{{\mathcal {U}}}

\def\P2{{P^{[2]}}}

\def\11{{\mbox{\boldmath $1$}}}

\def\eq{\begin{equation}}
\def\en{\end{equation}}

\def\eq#1\en{\begin{equation}#1\end{equation}}
\def\eqa#1\ena{\begin{eqnarray}#1\end{eqnarray}}

\newcommand{\OO}{\mathcal{O}}

\newcommand{\B}{{\mathbb B}}
\newcommand{\E}{{\mathbb E}}

\begin{document}

\title{The classifying topos of a topological bicategory}
\author{Igor Bakovi\'c \\Department of Mathematics\\Faculty of Natural Sciences and Mathematics\\ University of Split\\
Teslina 12/III, 21000 Split, Croatia\\ \\ and
\\ \\Branislav Jur\v co\\ Max Planck Institute for Mathematics\\Vivatsgasse 7,
53111 Bonn, Germany}
\date{}

\maketitle
\begin{abstract}
For any topological bicategory $\B$, the Duskin nerve $N\B$ of $\B$ is a simplicial space. We introduce the classifying topos ${\mathcal B}\B$
of $\B$ as the Deligne topos of sheaves $Sh(N\B)$ on the simplicial space $N\B$. It is shown that the category of geometric morphisms ${\rm
\Hom}(Sh(X), {\mathcal B}\B)$ from the topos of sheaves $Sh(X)$ on a topological space $X$ to the Deligne classifying topos is naturally
equivalent to the category of principal $\B$-bundles. As a simple consequence, the geometric realization $|N\B|$ of the nerve $N\B$ of a locally
contractible topological bicategory $\B$ is the classifying space of principal $\B$-bundles, giving a variant of the result of Baas, B\"okstedt
and Kro derived in the context of bicategorical K-theory. We also define classifying topoi of a topological bicategory $\B$ using sheaves on
other types of nerves of a bicategory given by Lack and Paoli, Simpson and Tamsamani by means of bisimplicial spaces, and we examine their
properties.
\end{abstract}
\thispagestyle{empty}

\footnotetext[1]{This research was in part supported by the Croatian Ministry of Science, Education and Sport, Project No. 098-0982930-2990.}

\footnotetext[2]{The first author acknowledges support from European Commission under the MRTN-CT-2006-035505 HEPTOOLS Marie Curie Research
Training Network}

\footnotetext[3]{Authors email addresses: ibakovic@gmail.com, branislav.jurco@googlemail.com}
\newpage
\section{Introduction}

In a recent paper by Baas, B\"okstedt and Kro \cite{BBK} it was shown that the geometric realization $|N\B|$ of the Duskin nerve $N\B$ \cite{D}
of a good topological bicategory $\B$ is the classifying space of charted $\B$-bundles. The bicategory is called good if its Duskin nerve $N\B$
is a good simplicial space, i.e. all degeneracy maps are closed cofibrations. Special cases of topological 2-groups and Lie 2-groups were
discussed in \cite{BS} and in \cite{J}, respectively.

The result of \cite{BBK} generalizes the well known fact that the geometric realization $|N\BC|$ of the nerve $N\BC$ of a locally contractible
topological category $\BC$ is the classifying space of principal $\BC$-bundles (on  CW complexes). This is very nicely described by Moerdijk in
\cite{M}. There, also the classifying topos ${\mathcal B}\BC$ of a topological category $\BC$ is described as the Deligne topos of sheaves
$Sh(N\BC)$ on the nerve $N\BC$ and it is shown that the category of geometric morphisms ${\rm \Hom}(Sh(X), {\mathcal B}\BC)$ from the topos of
sheaves $Sh(X)$ on a topological space $X$ to the Deligne topos is naturally equivalent to the category of principal $\BC$-bundles. As a simple
consequence, it is shown that the geometric realization $|N\BC|$ of the nerve $N\BC$ of a locally contractible topological category $\BC$ is the
classifying space of principal $\BC$-bundles.

One purpose of this note is to introduce the classifying topos ${\mathcal B}\B$ of a topological bicategory $\B$ as the topos of sheaves
$Sh(N\B)$ on the Duskin nerve $N\B$ of the bicategory $\B$, which is a simplicial space. The category of geometric morphisms ${\rm \Hom}(Sh(X),
{\mathcal B}\B)$ from the topos of sheaves $Sh(X)$ on a topological space $X$ to the classifying topos is naturally equivalent to the category
of (suitably defined) principal $\B$-bundles. As a simple consequence, the geometric realization $|N\B|$ of the nerve $N\B$ of a locally
contractible topological bicategory $\B$ is the classifying space of principal $\B$-bundles. Hence, we have a variant of the result of Baas,
B\"okstedt and Kro.

Another purpose of this note is to define classifying topoi of a topological bicategory $\B$ using sheaves on other types of nerves of the
bicategory $\B$, the nerves according to Lack \& Paoli \cite{LP} (or Simpson \cite{S} and Tamsamani \cite{T}), which can be viewed as bisimplicial
spaces. Again, the category of topos morphisms from the topos of sheaves $Sh(X)$ on a topological space $X$ to the corresponding classifying
topos is naturally equivalent to the respective category of (suitably defined) principal $\B$-bundles. As a simple consequence, the geometric
realization of any of these nerves of a locally contractible topological bicategory $\B$ is the classifying space of the respective principal
$\B$-bundles.

In Section 2, we recall some prerequisites from \cite{M} regarding sheaves on a simplicial space and augmented linear orders over topological
spaces. In Section 3, we recall, again from \cite{M}, the known facts about classifying spaces and topoi of topological categories (and the
corresponding principal bundles). We describe a generalization to the case of bicategories, based on the Duskin nerve, in Section 4. Further
preliminaries needed for the subsequent discussion of alternative definitions of classifying spaces and topoi of bicategories are given in
Section 5. Finally, in section 6, we describe a modification of the classifying topos of a topological bicategory (and the corresponding
principal bundles) based on alternative definitions of the nerves according to Lack \& Paoli, Simpson and Tamsamani.

This article is meant to be the first one in the sequel within a program, initiated by the authors, of classifying topoi of higher order
structures in topology. It is a vast generalization of the program initiated by Moerdijk in \cite{M} on the relation between classifying spaces
and classifying topoi. Moerdijk's lecture notes arose out of an important question: What does the classifying space of a small category
classify? In the article titled by the same question \cite{W}, Weiss proved the classifying property of the classifying space for slightly
different geometric objects then those of Moerdijk, showing that the answer may not be unique.

Therefore, this article may be seen as an (one possible) answer to the following question: What does the classifying space of a topological
bicategory classify? Bicategories are the weakest possible generalization of ordinary categories to the immediate next level of dimension. Like
categories, bicategories do have a genuine simplicial set associated with them, their Duskin nerve \cite{D}. Based on unpublished work of
Roberts on the characterization of the nerve of a strict n-category, Street postulated in \cite{St} an equivalence between the category of {\it
strict $\omega$-categories} and a category of certain types of simplicial sets which are called {\it complicial sets}. The {\it Street-Roberts
conjecture} was proved by Verity in \cite{V1}, and in his subsequent papers \cite{V2} and \cite{V3} he gave a characterization of {\it weak
$\omega$-categories}. Under this characterization, one should be able to capture classifying spaces and topoi of bicategories and other higher
dimensional categories, at least in so far as these concepts have found satisfactory definitions. Following such reasoning, we may define the
classifying space of a weak $\omega$-category as a geometric realization of the complicial set which is its nerve, and the classifying topos of
a weak $\omega$-category as a topos of sheaves on that complicial set.

It would be interested to compare this approach to classifying spaces of weak $\omega$-categories with classifying spaces of crossed complexes
defined by Brown and Higgins in \cite{BH2}, since there is a well known equivalence between {\it strict $\omega$-groupoids} and {\it crossed
complexes} proved in \cite{BH1} by the same authors. In particular, it would be interesting to see whether the methods we developed would allow
to define a classifying space of a weak $\omega$-category by taking a fundamental crossed complex of its coherent simplicial nerve.

However, this article is not so cosmological in its scope, and its main contribution is to put together some established results on classifying
spaces and classifying topoi in a new way, with consequences for the theory of bicategories. Since we are following Moerdijk's approach to
classifying spaces and classifying topoi, we will omit all proofs, which can be found in Moerdijk's lecture notes.

{\bf Acknowledgments} We would like to thank Ronnie Brown for useful comments on the relevance of this work to general notions of classifying
spaces of crossed complexes. We would also like to express our thanks to the referee of "Homology, Homotopy and Applications" whose careful
reading resulted in comments and suggestions which have improved the structure and the content of this article. The first author would like to
thank for a hospitality of Max Planck Institute for Physics in Munich, where the part of this research was made, and to Max Planck Institute for
Mathematics in Bonn, where he was supported by the Croatian-German bilateral DAAD program "Homological algebra in geometry and physics".

\section{Simplicial spaces and linear orders over topological spaces}
In this section, we recall some prerequisites regarding sheaves on a simplicial space and augmented linear orders over topological spaces.
Almost all the definitions and theorems are taken verbatim from \cite{M}, where proofs of all statements of this section can be found.

\subsection{Topological spaces}

Let us recall that a closed set in a (topological) space $X$ is {\it irreducible} if it can not be written as a union of two smaller closed
sets. The space $X$ is {\it sober} if every irreducible set is the closure $\overline{\{x\}}$ of the one point set $\{x\}$ of a unique $x\in X$.
Every Hausdorff space is sober. In this note all spaces will be sober by assumption.

A space $X$ is {\it locally equiconnected} (LEC) if the diagonal map $X\to X \times X$ is a closed cofibration. For example, CW-complexes are
LEC.

A space $X$ is {\it locally contractible} if it has a basis of contractible sets. Examples of locally contractible spaces are locally
equiconnected spaces and in particular CW complexes. For a locally contractible space the \'etale homotopy groups $\pi_n(Sh(X), x_0 )$ are
naturally isomorphic to the ordinary homotopy groups $\pi_n(X, x_0 )$ for each $n$.

\subsection{Sheaves as \'etale spaces}

Throughout this article, we will consider sheaves as sheaves of cross-sections of \'etale spaces. Recall that a bundle $p\colon E \to X$ over
$X$ is said to be {\it \'etale space} over $X$ if for each $e \in E$ there exists an open set $V \subset E$, with $e \in V$, such that $p(V)
\subset X$ is open in $X$ and the restriction $p|_V\colon V \to p(V)$ over $V$ is a homeomorphism. There is a well known equivalence
\[\xymatrix{Etale(X) \ar@<0.5ex>[r]^-{\Gamma} & Sh(X) \ar@<0.5ex>[l]^-{\Lambda}}\]
where $\Gamma \colon Etale(X) \to Sh(X)$ is a functor which assigns to each \'etale space $p\colon E \to X$ over $X$ the sheaf of all
cross-sections of $E$. The functor $\Lambda \colon Sh(X) \to Etale(X)$ assigns to each sheaf $S$ the \'etale space of {\it germes} of $S$, where
the germ at the point $x \in X$ is an equivalence class $germ_x s$ represented by $s \in S(U)$ under the equivalence relation, which relates two
elements $s \in S(U)$ and $t \in S(V)$ if there is some open set $W \subset U \cap V$ such that $x \in W$ and $s|_{W}=t|_{W}$. The stalk of the
sheaf $S$ at the point $x \in X$ is the set $S_x=\{germ_x s \colon s \in S(U), x \in U\}$ of all germs at $x$, which is formally a filtered
colimit
\[ S_x = \underset{x \in U}{\underrightarrow{lim}} S(U) \] of the restriction $S^{(x)} \colon \OO_x(X)^{op} \to Set$ of the sheaf $S$ to the filtered
category $\OO_x(X)^{op}$ of open neighborhoods of the point $x \in X$. Then $\Lambda S$ is an \'etale space $p \colon \coprod_{x \in X} S_x \to
X$ whose sheaf of cross sections is canonically isomorphic to $S$. Therefore, we will use simultaneously terms sheaves and \'etale spaces in the
rest of this article.

\subsection{Topoi} I n the following, a  {\it topos}  will always mean a Grothendieck topos. $Sh(X)$ will denote topos of sheaves on a (topological)
space $X$. A sober space $X$ can be recovered from the topos $Sh(X)$, which is the faithful image of the space $X$ in the world of topoi.

Further, ${\rm Hom}(Sh(X), Sh(Y))$ will denote the category of geometric morphisms from $Sh(X)$ to $Sh(Y)$. We will use the same notation ${\rm
Hom}({\mathcal F}, {\mathcal E})$ also in the more general case of any two topoi ${\mathcal F}$ and ${\mathcal E}$. By definition a {\it
geometric morphism} $f \in {\rm Hom}({\mathcal F}, {\mathcal E})$ is a pair of functors $f^*\colon{\mathcal E}\to {\mathcal F}$ and
$f_*\colon{\mathcal F}\to {\mathcal E}$,  $f^*$ being left adjoint to $f_*$ and also $f^*$ being left exact, i.e. preserving finite limits.

Let us recall that a geometric morphism $f\colon {\mathcal F}\to {\mathcal E}$ between locally connected topoi is a {\it weak homotopy
equivalence} if it induces an isomorphisms on \'etale homotopy (pro)groups $\pi_0({\mathcal F})\cong \pi_0({\mathcal E})$ and $\pi_n({\mathcal
F, p})\cong \pi_n({\mathcal E, fq})$, for $n\geq 1$ for any base point $q \in {\mathcal F}$.

For the collection of homotopy classes of geometric morphism from ${\mathcal F}$ to ${\mathcal E}$ the usual notation $[{\mathcal F}, {\mathcal
E}]$ will be used.
\subsection{The singular functor} \label{sing} The following construction
of a singular functor is taken from \cite{Kelly}, where Kelly described it in the context of enriched $\mathcal V$-categories for any symmetric
monoidal closed category $\mathcal V$, which is complete and cocomplete. Let
\[\label{singemb} F\colon {\BA} \to
{\E}
\]
be a functor from the small category $\BA$. The {\it singular functor} of $F$ is the functor
\[\label{singfun} {\E}(F,1) \colon
{\E} \to [{\BA}^{op},{\mathcal V}]
\] which is
obtained as the composite of the Yoneda embedding
\[\label{Yonemb}
{\rm Yon}\colon{\E} \to [{\E}^{op},{\mathcal V}]
\]
followed by the functor $[F^{op},{\mathcal V}] \colon[{\E}^{op},{\mathcal V}] \to [{\BA}^{op},{\mathcal V}]$ given by restriction along a
functor $F$. More explicitly, the singular functor ${\E}(F,1)$ sends any object $E$ in $\E$ to the functor
\[\label{singfunob}
{\E}(F(-),E) \colon {\BA}^{op} \to {\mathcal V}
\]
which takes an object $A$ in ${\BA}$ to the hom-object ${\E}(F(A),E)$ in ${\mathcal V}$. If the category $\E$ is cocomplete, then the singular
functor has a left adjoint
\[\label{tensprod} L \colon [{\BA}^{op},{\mathcal V}] \to \E \]
defined for each presheaf $P \colon {\BA}^{op} \to {\mathcal V}$ as the colimit
\[\label{tenscolimit}\xymatrix{ L(P)={\underrightarrow{lim}}(\int_{\BA} P \ar[r]^-{\pi_{P}} & \BA \ar[r]^{F} & \E)}\]
where $\int_{\BA} P$ is the so called Grothendieck construction \cite{MM} on a presheaf $P \colon {\BA}^{op} \to {\mathcal V}$.
\subsection{Grothendieck nerve as a singular functor}\label{Groth sing}
Each ordinal $[n]=\{0 < 1 < \ldots <n \} $ can be seen as a category with objects $0,1,\ldots n,$ and a unique arrow $i \to j$ for each $0 \leq
i \leq j \leq n$. Also, any monotone map between two ordinals may be seen as a functor. In this way, $\Delta$ becomes a full subcategory of
$\mathbf{\rm Cat_1}$ with a fully faithful inclusion functor
\[\label{Deltacat1} J\colon \Delta \to \mathbf{\rm
Cat_1}
\]
For any small category $\B$, the composite of the Yoneda embedding ${\rm Yon}\colon\B \to [\B^{\rm op},\mathbf{\rm Set}]$ followed by the
restriction functor $[\B^{\rm op},\mathbf{\rm Set}] \to [\Delta^{\rm op},\mathbf{\rm Set}]$ along $J$ gives a singular functor of $J$. In more
details, the singular functor of $J$ defines the Grothendieck nerve functor
\[\label{Gronerve} N\colon \mathbf{\rm
Cat_1}\to [\Delta^{\rm op},\mathbf{\rm Set}]
\]
which sends any category $\BC$ to the simplicial set $N\BC$ which is the {\it nerve} of $\BC$ whose n-simplices are defined by the set
\[\label{nsimpGronerve} N\BC_{n}=[J([n]),\BC]
\]
where the right side denotes the set of functors from an ordinal $[n]$ to the category $\BC$. The nerve functor is fully faithful, which means
that the simplicial skeletal category $\Delta$ is an adequate subcategory of the category $\mathbf{\rm Cat_1}$ in the sense of Isbell
\cite{Is1}, \cite{Is2}. We also say that the corresponding embedding is dense, in the sense of Kelly \cite{Kelly}.

\subsection{Simplicial spaces}\label{simpsp} Let $\Delta$ be the simplical model
category having as objects nonempty finite sets (ordinals) $[n]=\{0,1,\ldots,n\}$, for $n\geq 0$, and as arrows order-preserving functions
$\alpha\colon [n]\to [m]$. A {\it simplicial space (set)} is a contravariant functor from $\Delta$ into the category of spaces (sets). Its value
at $[n]$ is denoted $Y_n$ and its action on arrow $\alpha\colon [n]\to [m]$ as $Y(\alpha)\colon Y_m\to Y_n$. A simplicial space $Y$ is called
{\it locally contractible} if each $Y_n$ has a basis of contractible sets.

For a simplicial space $Y$ the {\it geometric realization} $|Y|$ will always mean the thickened (fat) geometric realization. This is defined as
a topological space obtained from the disjoint sum $\sum_{n\geq 0} X_n \times \Delta^n$ by the the equivalence relations
\[ (\alpha^*(x), t) \sim (x, \alpha(t))\]
for all injective (order-preserving) arrows $ \alpha\colon [n]\to [m] \in \Delta$, any $x \in X_m$ and any $t\in \Delta^n$, where $\Delta^n$ is
the standard topological $n$-simplex. If all degeneracies are closed cofibrations, i.e. the simplicial space is a {\it good} simplicial space,
this geometric realization is homotopy equivalent to the geometric realization of the underlying simplicial set of $Y$, which is defined as
above but allowing for all arrows in $\Delta$. In particular, $Y$ is good if all spaces $Y_n$ are locally equiconnected \cite{BBK}. Geometric
realization of a locally contractible simplicial space is a locally contractible space.

\begin{definition}\label{Sh(Y)}A {\bf sheaf $S$ on a simplicial space $Y$}
is defined to be a system of sheaves $S^n$ on $Y_n$, for $n\geq 0$, together with sheaf maps $S(\alpha)\colon Y(\alpha)^* S^n\to S^m$ for each
$\alpha \colon [n] \to [m]$. These maps are required to satisfy the following functoriality conditions:

\begin{itemize}
\item [i)] (normalization) $S(\rm{id_{[n]}})= \rm{id_{S_{n}}}$, and

\item [ii)] for any $\alpha \colon [n]\to [m]$, $\beta\colon [m]\to [k]$ the following diagram
\[\xymatrix@!=4pc{Y(\beta)^*Y(\alpha)^* S^n \ar[d]_{\cong}
\ar[rr]^-{Y(\beta)^*S(\alpha)} && Y(\beta)^* S^m \ar[d]^{S(\beta)}\\
Y(\beta\alpha)^* S^n \ar[rr]_{S(\beta\alpha)} && S^k}\]
\end{itemize}
is commutative.

\noindent A morphism $f \colon S\to T$ of sheaves on $Y$ consists of maps $f_n \colon S^n\to T^n$ of sheaves on $Y_n$ for each $n\geq 0$, which
are compatible with the structure maps $S(\alpha )$ and $T(\alpha)$. This defines the category $Sh(Y)$ of sheaves on the simplicial space $Y$.
\end{definition}
\begin{proposition} \label{Deligne1} The category $Sh(Y)$ of sheaves
on a simplicial space is a topos.
\end{proposition}
\begin{theorem} \label{weak-equiv1} For any simplicial space $Y$ the topoi
$Sh(Y)$ and $Sh(|Y|)$ have the same weak homotopy type.
\end{theorem}
\begin{definition} A {\bf linear order over a topological space $X$} is a
sheaf $p \colon L \to X$ on $X$ together with a subsheaf $O \subseteq L\times_X L$ such that for each point $x \in X$ the stalk $L_x$ is
nonempty and linearly ordered by the relation
\[ y\leq z \hskip1cm  \mbox{iff} \hskip1cm (y,z)\in O_x,\]
for $y,z\in L_x$. A mapping $L \to L'$ between two linear orders over $X$ is a mapping of sheaves restricting for each $x\in X$ to an order
preserving map of stalks $L_x \to L'_x$. This defines a category of linear orders on $X$.
\end{definition}
\begin{example} An {\it open ordered covering} ${\mathcal U}=\{U_i\}_{i \in I}$ of a topological space $X$, is a covering indexed over a
partially ordered set $I$, which restricts to a total ordering on every finite subset $\{i_0,\ldots,i_n\}$ of $I$ whenever the finite
intersection $U_{i_0,\ldots,i_n}=U_{i_0} \cap \ldots \cap U_{i_n}$ is nonempty. When a sheaf $p \colon L \to X$ is given by the projection $p
\colon \coprod_{i \in I} U_i \to X$ from the disjoint union of open sets in the open ordered covering ${\mathcal U}$ the subsheaf $p^{[2]}
\colon L \times_X L \to X$ is given by the induced projection $p^{[2]} \colon \coprod_{i,j \in I} U_{ij} \to X$ from the family $\{U_{ij}\}_{i,j
\in I}$ of double intersections of open sets ${\mathcal U}$. The family of inclusions $i_{ij} \colon U_{ij} \hookrightarrow \coprod_{i,j \in I}
U_{ij}$, for each $U_{ij} \neq \emptyset$ such that $i < j$, defines a subsheaf $O = \coprod_{i<j} U_{ij}$ of $L = \coprod_{i,j \in I} U_{ij}$
whose stalks $O_x$ are linearly ordered for any $x \in X$. Therefore, open ordered coverings used by Baas, B\"okstedt and Kro in \cite{BBK} are
examples of linear orders over $X$.
\end{example}
\begin{remark}A linear order
$L$ over $X$ defines an obvious topological category with $L$ as the space of objects and the order subsheaf $O \subseteq L\times_X L$ as the
space of arrows. Hence, we can speak of a nerve $NL$ of the linear order $L$. This nerve is obviously a simplicial sheaf on $X$ (a simplicial
space with \'etale maps into $X$).
\end{remark}
Recall that any open covering of a topological space $X$ can be assembled into a simplicial sheaf over $X$ with distinguished properties.
Therefore, by the construction in the Example 2.5 and following the Remark 2.6 we may regard linear orders as generalizations of coverings of
topological spaces.
\begin{definition}\label{aug}For any space $X$ and any simplicial space
$Y$ write $Lin(X,Y)$ for the category of linear orders $(L,\rm{aug})$ over $X$ equipped with a simplicial map ({\bf augmentation})
${\rm{aug}}\colon NL \to Y$ from the nerve of $L$ to $Y$. A morphism $(L,{\rm{aug})}\to (L',\rm{aug}')$ in $Lin(X,Y)$ are maps of linear orders
$L\to L'$ such that the induced map $NL\to NL'$ on the nerves respects the augmentations.
\end{definition}
If we regard linear orders as generalizations of coverings of topological spaces, then augmentations of linear orders may be seen as cocyles on
such coverings.
\begin{example} Let $N\BC$ be the nerve of a topological category $\BC$. An augmentation ${\rm{aug}}\colon NL \to N\BC$ of
a linear order $L$ defined by an open ordered covering ${\mathcal U}=\{U_i\}_{i \in I}$ of a topological space $X$ as in Example 2.5 is a \v
Cech cocyle on the covering ${\mathcal U}$ with values in the category $\BC$.
\end{example}
\begin{definition}\label{conc}
We call two objects $E_0, E_1 \in Lin(X,Y)$ {\bf concordant} if there exists an $ E\in Lin(X\times [0,1],Y)$ such that we have $E_0\cong
i^*_0(E)$ and $E_0\cong i^*_1(E)$ under the obvious inclusions $i_0, i_1\colon X \hookrightarrow X\times [0,1]$. $Lin_c(X,Y)$ will denote the
collection of concordance classes of objects from $Lin(X,Y)$.
\end{definition}
\begin{theorem}\label{equiv1}Let $Y$ be a simplicial space. For any space
$X$ there is a natural equivalence of categories
\[{\rm Hom}(Sh(X), Sh(Y))\simeq Lin(X,Y).\] On homotopy classes of topos morphisms we
have the natural bijection
\[[Sh(X), Sh(Y)]\cong Lin_c(X,Y).\]
\end{theorem}
\begin{corollary}\label{classifi} Let $X$ be a CW-complex and $Y$ be a locally contractible
simplicial space. There is a natural bijection between homotopy classes of maps $[X,|Y|]$ and concordance classes $Lin_c(X,Y)$.
\end{corollary}
\begin{remark} If in addition the simplicial space $Y$ is a good one then the
above is true also if we use, instead of its thickened geometric realization, the geometric realization  of the underlying simplicial set of
$Y$. In particular, it doesn't matter which geometric realization we use if each of $Y_n$ is LEC or a CW-complex.
\end{remark}
\section{Classifying spaces and classifying topoi of topological categories}
In this section we specify the known results described in Section 1 to the case when the simplicial space $Y$ is the nerve of a topological
category $\BC$. The reader who is interested in more details is referred to \cite{M}, which we again follow almost verbatim.
\begin{definition} \label{topTop} Let $\BC$
be a topological category. The {\bf classifying topos} ${\mathcal B}\BC$ of a topological category is defined as the topos $Sh(N\BC)$.
\end{definition}
\begin{definition} The {\bf classifying space} $B\BC$ of a topological category
$\BC$ is the geometric realization $|N\BC|$ of its nerve $N\BC$.
\end{definition}
\vskip0.5cm With these definitions we have the following corollary of Theorem \ref{weak-equiv1}.
\begin{corollary}\label{weak-equiv2}For any topological category $\BC$ the topos
of sheaves $Sh(B\BC)$ on the classifying space $B\BC$ has the same weak homotopy type as the classifying topos ${\mathcal B} \BC$.
\end{corollary}
\begin{definition} For any topological category $\BC$
write $Lin(X,\BC)$ for the category of linear orders over $X$ equipped with an augmentation  $NL\to N\BC$. An object $E$ of this category will
be called a {\bf principal} $\BC$-{\bf bundle}. We call two principal $\BC$-bundles $E_0$ and $E_1$ on $X$ {\bf concordant} if there exists a
principal $\BC$-bundle on $X\times [0,1]$ such that we have isomorphisms $E_0\cong i^*_0(E)$ and $E_0\cong i^*_1(E)$ under the obvious
inclusions $i_0, i_1\colon X \hookrightarrow X\times [0,1]$.
\end{definition}
\begin{remark} The nerve construction leads to a bijection between
principal $\BC$-bundles and linear orders $L$ equipped with a continuous functor $L\to \BC$.
\end{remark}
\vskip0.5cm The fact that the classifying topos ${\mathcal B}\BC$ classifies principal $\BC$-bundles follows now immediately from Theorem
\ref{equiv1}.
\begin{theorem} \label{equiv2}For a topological category $\BC$ and a
topological space $X$ there is a natural equivalence of categories
\[{\Hom(Sh(X),{\mathcal B} \BC})\simeq Lin(X,\BC).\] On homotopy classes of topos morphisms we
have the natural bijection
\[[Sh(X), {\mathcal B} \BC]\cong Lin_c(X,\BC).\]
\end{theorem}
\vskip0.5cm Similarly, the fact that the classification space $B\BC$ classifies principal $\BC$-bundles now follows from Corollary
\ref{classifi}.
\begin{definition} We say that a topological category $\BC$ is {\bf locally contractible} if its space of objects $C_0$ and its space of
arrows $C_1$ are locally contractible. A topological category $\BC$ is a {\bf good topological category}, if its nerve $N\BC$ is a good
simplicial space.
\end{definition}
\begin{corollary}\label{classifi2}
For a locally contractible category $\BC$ and a CW-complex $X$ there is a natural bijection
\[[X, B\BC] \cong Lin_c(X,\BC).\]
\end{corollary}
\begin{remark} If, in addition, the topological category $\BC$ is a good one then the above is true also if we use, instead of the thickened geometric
realization of the nerve, the geometric realization  of the underlying simplicial set. In particular, it doesn't matter which geometric
realization we use if all $N\BC_n$ are LEC.
\end{remark}

\section{Classifying spaces and classifying topoi of topological bicategoies I}
In this section we specify the known results described in Section 1 to the case when the simplicial space $Y$ is the nerve of a topological
bicategory $\B$.
\subsection{Duskin nerve as a singular functor}\label{Duskin sing}
The Duskin nerve \cite{D} can also be obtained as a singular functor when we take ${\mathcal V}=\mathbf{\rm Set}$. Every category (in particular
the category defined above by the ordinal $[n]$) can be seen as a locally discrete bicategory (the only 2-cells are identities) which gives a
fully faithful inclusion
\[\label{Deltabicat1} H \colon \Delta \to \mathbf{\rm Bicat_1}
\]
where $\mathbf{\rm Bicat_1}$ denotes the category of bicategories and normal lax functors, or normal morphisms of bicategories defined by
B\'enabou in \cite{Be}. The singular functor of the inclusion $H$ is the Duskin nerve functor
\[\label{Dusnerve} N\colon \mathbf{\rm
Bicat_1} \to [\Delta^{\rm op},\mathbf{\rm Set}]
\]
which is fully faithful and sends a (small) bicategory $\B$ to its nerve $N\B$ which is a simplicial set whose n-simplices are defined by the
set
\[\label{nsimpDusnerve} N\B_{n}=[H([n]),\B]
\]
where right side is a set of normal lax functors from an ordinal $[n]$ to the bicategory $\B$.

\begin{definition} For an ordinal $[n]$ and a bicategory $\B$ a
{\bf normal lax functor} $(B,f,\beta)\colon [n]\to \B$ consists of the following data in $\B$:
\begin{itemize}
\item [(i)] an object $B_i$ for each $i\in [n]$,

\item [(ii)] a morphism $f_{ij}\colon B_i \to B_j$ for each $i, j \in [n]$ with
$i \leq j$,

\item [(iii)] a 2-cell $\beta_{ijk} \colon f_{ij} \circ f_{jk}
\Rightarrow f_{ik}$ for each $i, j, k\in [n]$ with $i \leq j \leq k$
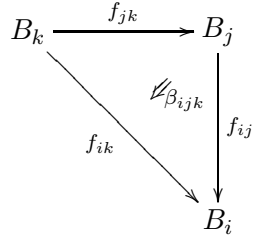
\[\xymatrix@!=4pc{B_{k} \drtwocell<\omit>{<-4>\,\,\,\,\,\,\,\beta_{ijk}}
\ar[r]^-{f_{jk}} \ar[dr]_-{f_{ik}} & B_{j} \ar[d]^-{f_{ij}}\\
& B_{i}}\]
\end{itemize}
such that the following conditions are satisfied:
\begin{itemize}
\item ({\it normalization}) for any $i \in [n]$ we have $f_{ii}=i_{B_i} \colon B_i \to
B_i$ and for any $i, j \in [n]$ such that $i \leq j$ the corresponding 2-cells $\beta_{iij} \colon f_{ii} \circ f_{ij} \Rightarrow f_{ij}$  and
$\beta_{ijj} \colon f_{ij} \circ f_{jj} \Rightarrow f_{ij}$ are given by the two 2-simplices
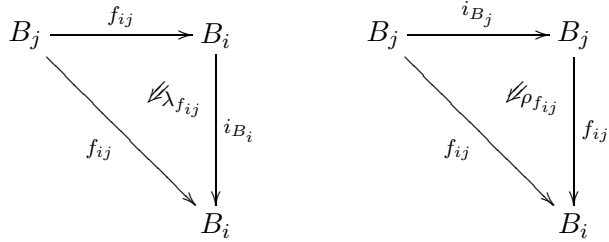
\[\xymatrix@!=4pc{B_{j}
\drtwocell<\omit>{<-4>\,\,\,\,\,\,\,\lambda_{f_{ij}}}
\ar[r]^-{f_{ij}} \ar[dr]_-{f_{ij}} & B_{i} \ar[d]^-{i_{B_{i}}}\\
& B_{i}}  \hspace{3pc} \xymatrix@!=4pc{B_{j} \drtwocell<\omit>{<-4>\,\,\,\,\,\,\,\rho_{f_{ij}}}
\ar[r]^-{i_{B_{j}}} \ar[dr]_-{f_{ij}} & B_{j} \ar[d]^-{f_{ij}}\\
& B_{i}}\] where $\rho_{f_{ij}} \colon f_{ij} \circ i_{p_{j}} \Rightarrow f_{ij}$ and $\lambda_{f_{ij}} \colon i_{p_{i}} \circ f_{ij}
\Rightarrow f_{ij}$ are the components of the right and left identity natural isomorphisms in $\B$.

\item ({\it coherence condition}) for each $i, j, k, l\in [n]$ such that $i \leq j \leq k \leq
l$ the following tetrahedron
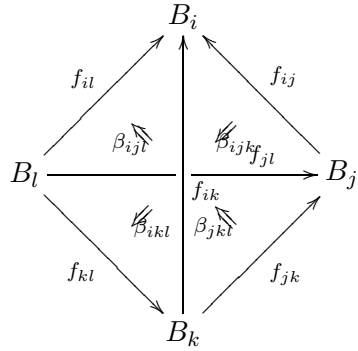
\[\xymatrix@!=3pc{& B_{i} \drtwocell<\omit>{<4>\beta_{ijk}} \dltwocell<\omit>{<-4>\beta_{ijl}} &\\
B_{l} \drtwocell<\omit>{<-4>\beta_{ikl}} \ar[ur]^{f_{il}} \ar'[r][rr]^{f_{jl}} \ar[dr]_-{f_{kl}}  &&
B_{j} \dltwocell<\omit>{<4>\beta_{jkl}} \ar[ul]_-{f_{ij}} \\
& B_{k} \ar[uu]_(0.45){f_{ik}} \ar[ur]_-{f_{jk}} &}\] is commutative. This means that we have the identity of 2-cells in the bicategory $\B$:
\[\label{Dusnervecoh} \beta_{ikl}(\beta_{ijk}
\circ f_{kl})=\beta_{ijl} (f_{ij} \circ \beta_{jkl}) \alpha_{ijkl}.
\]
\end{itemize}
\end{definition}
\begin{remark} Simplicial sets that are isomorphic to a nerve of
a bicategory have been characterized in \cite{D} and \cite{G}. Simplicial sets that are isomorphic to a nerve of a bicategory form a full
subcategory of the category of simplicial sets. This category is equivalent to the category ${\rm Bicat}_1$ of bicategories with lax normal
functors. Let us recall that a lax functor $(F, \phi)$ is normal if $F(id_x)=id_{Fx}$ and $\phi_x\colon id_{Fx}\implies F(id_x)$ is the identity
2-cell and oplax means that all the structure maps go in opposite direction. This equivalence holds also in the topological setting.
\end{remark}
\begin{definition} Let $\B$
be a topological bicategory. The {\bf classifying topos} ${\mathcal B}\B$ of the topological bicategory $\B$ is defined as the topos $Sh(N\B)$.
\end{definition}
\begin{definition} The {\bf classifying space} $B\B$ of a topological bicategory
$\B$ is the geometric realization $|N\B|$ of its nerve $N\B$.
\end{definition}
\vskip0.5cm With these definitions we have the following corollary of Theorem \ref{weak-equiv1}.

\begin{corollary} For any topological bicategory $\B$ the topos of
sheaves $Sh(B\B)$ on the classifying space $B\B$ has the same weak homotopy type as the classifying topos ${\mathcal B}\B$.
\end{corollary}
\begin{definition} For a topological bicategory $\B$
write $Lin(X,\B)$ for the category of linear orders over $X$ equipped with an augmentation ${\rm{aug}}\colon NL \to N\B$. An object $E$ of this
category will be called a {\bf Duskin principal} $\B$-{\bf bundle}. We call two Duskin principal $\B$-bundles $E_0$ and $E_1$ on $X$ {\bf
concordant} if there exists a Duskin principal $\B$-bundle on $X\times [0,1]$ such that we have the equivalences $E_0\simeq i^*_0(E)$ and
$E_0\simeq i^*_1(E)$ under the obvious inclusions $i_0, i_1\colon X \hookrightarrow X\times [0,1]$.
\end{definition}
\begin{remark} We can consider a linear order $L$ as a locally trivial bicategory
(with only trivial 2-morphisms). In this case the Duskin nerve of $L$ coincides with the ordinary nerve of $L$ which justifies the same notation
$NL$ for both nerves.
\end{remark}
\begin{remark} By the above
remark an augmentation $NL \to N\B$ is the same, by the nerve construction, as a continuous normal lax functor $L\to \B$.
\end{remark}
\vskip0.5cm Similarly to Theorem \ref{equiv2} we have from Theorem \ref{equiv1} the following "classifying" property of the classifying 1-topos
${\mathcal B} \B.$
\begin{theorem}\label{equivDus}For a topological bicategory $\B$ and a
topological space $X$ there is a natural equivalence of categories
\[{\Hom (Sh(X),{\mathcal B} \B})\simeq Lin(X,\B).\] On homotopy classes of topos morphisms we
have the natural bijection
\[[Sh(X), {\mathcal B} \B]\cong Lin_c(X,\B).\]
\end{theorem}
\begin{definition} We say that a topological bicategory $\B$ is {\bf locally contractible} $\B$ if its space of objects $B_0$, its space of
1-arrows $B_1$ and its space of 2-arrows $B_2$ are locally contractible. A topological bicategory $\B$ is a {\bf good topological bicategory},
if its nerve $N\B$ is a good simplicial space.
\end{definition}
\vskip0.5cm The "classification" property of the classifying space $B\B$ now follows as a corollary from Corollary \ref{classifi}.
\begin{corollary} For a locally contractible bicategory $\B$ and a CW-complex $X$ there is a
natural bijection
\[[X, B\B] \cong Lin_c(X,\B).\]
\end{corollary}
\begin{remark}
If, in addition, the topological bicategory $\B$ is a good one then the above is true also if we use, instead of the thickened geometric
realization of the nerve, the geometric realization  of the underlying simplicial set.  In particular,  it doesn't matter which geometric
realization we use if all $N\B_n$ are LEC. The case of a good topological bicategory, as well as the sufficient conditions for a bicategory
being a good one, are discussed in \cite{BBK}. Those conditions actually guarantee that all $N\B_n$ are LEC. Thus, our corollary above gives a
slight generalization of the result of Baas, B\"okstedt and Kro.
\end{remark}

\section{Principal bundles under a category}

Before introducing an alternative notion of a classifying topos of a bicategory in the next section, we will introduce some additional background
material. Everything up to including remark \ref{si} is taken almost verbatim from \cite{M} where the interested reader can find the missing
proofs (as well as more details). The Definition \ref{def}, Theorem \ref{equiv10} and Corollary \ref{classifi13} might be new. To make our
discussion more complete we start with definitions of the classifying topoi in cases of a small and s-\'etale category. We also recall the
definition of a principal $\BC$-bundles in these cases.
\begin{proposition}The category of all presheaves on
a small category $\BC$ is a topos.
\end{proposition}
\begin{definition} \label{smallTop} The topos ${{\mathcal B}\BC}$ of presheaves on
a small category $\BC$ is called the {\bf classifying topos} of $\BC$. \vskip0.5cm
\end{definition}
\begin{remark}
At this point, the reader may wonder how the above definition \ref{smallTop} is related to the definition of the classifying topos of a
topological category given in Definition \ref{topTop}. We will address this question later in \ref{equivDefTop}, after we introduce further
relevant material.
\end{remark}

\newpage
\begin{definition}\label{principal} For a small category $\BC$ and a space $X$,
a {\bf $\BC$-bundle over $X$} is a covariant functor $E\colon \BC\to Sh(X)$. Such a $\BC$-bundle is called a {\bf principal} (flat, filtering)
if for each point $x\in X$ the following conditions - non-emptiness, transitivity and freeness - are satisfied for the stalks $E(c)_x$ for
objects $c\in \BC$ :

(i) There is at least one object $c$ in $\BC$ for which the stalk $E(c)_x$ is non-empty.

(ii) For any two points $y\in E(c)_x$ and $z\in E(d)_x$, there are arrows $\alpha\colon b\to c$ and $\beta\colon b\to d$ from some object $b$ of
$\BC$, and an object $w\in E(b)_x$ such that $\alpha w=y$ and $ \beta w=z$.

(iii) For any two parallel arrows $\alpha, \beta\colon c \to d$ and any $y\in E(c)_x$ for which $\alpha y = \beta y$ there is an arrow
$\gamma\colon b\to c$ and a point $z\in E(b)_x$ such that $\alpha \gamma = \beta \gamma$ and $\gamma z = y$.

A map between two principal $\BC$-bundles is a natural transformation between the corresponding functors. The category of principal
$\BC$-bundles will be denoted as ${\rm Prin}(X,\BC)$.
\end{definition}

\begin{examples} The following well known notions are examples of principal $\BC$-bundles:

(i) (principal group bundles) Any group $G$ can be seen as a groupoid (and therefore a category) with only one object. In this way, the above
definition of a principal $\BC$-bundle becomes the usual one where a {\it principal left $G$-bundle} over $X$ is an \'etale space $p \colon P
\to X$ with a with a fibre-preserving left action $a \colon G \times P \to P$ of $G$ on $P$ for which the induced map $(a,pr_2) \colon G \times
P \to P \times P$ is a homeomorphism.

(ii) (principal monoid bundles) Any monoid $M$ can be seen as a category with only one object. If every morphism in a such category is a
monomorphism, then the monoid $M$ is said to have left cancelation if $mk=ml$ implies $k=l$ for any $k,l,m \in M$. Segal used such a monoid $M$
in order to introduce a {\it right principal monoid bundle} in \cite{Se} as an \'etale space $p \colon P \to X$ over $X$ with a fibre-preserving
right action of $M$ on $P$, such that each fibre $P_x$ is a principal $M$-set. A {\it right principal $M$-set} $S$ is a set with a right action
of $M$ which is free in the sense that $s m_1=s m_2$ for any $m_1,m_2 \in M$ and $s \in S$, and transitive in the sense that for any $s_1,s_2
\in M$ there exist $m_1,m_2 \in M$ and $s \in S$ such that $s_1=s m_1$ and $s_2=s m_2$. Although Segal used right action of a monoid with left
cancelation, it is obvious that when $\BC$ is a monoid $M$ with right cancelation, the above definition of a left principal $\BC$-bundle becomes
left principal monoid bundle.

(iii) (principal poset bundles) Any partially ordered set $P$ may be seen as a category with exactly one morphism $i \to j$ if and only if $i
\leq j$. A principal $P$-bundle over a topological space $X$ is a covering $\U=\{U_{i}\}_{i \in P}$ of $X$ such that when $i \leq j$ then $U_{i}
\subseteq U_{j}$ and which is {\bf locally directed} in the sense that any $U_{ij}=U_{i} \cap U_{j}$ is covered by the family $\U_{ij}=\{U_{k}
\colon k \leq i \, \wedge \, k \leq j\}$.

(iv) (principal simplicial sets) Any linear order over a topological space defines a topological category and therefore a simplicial space via
its nerve as in the Remark 2.6. One can see that a simplicial set $S \colon \Delta^{op} \to Set$ is a principal $\Delta^{op}$-bundle if and only
if is the nerve of a (uniquely determined) non-empty linear order.
\end{examples}

\newpage
\begin{definition} A $\BC${\bf-sheaf} is an \'etale space $p \colon S \to C_0$ equipped with a continuous right $\BC$ action
\[\alpha \colon S\times_{C_0} C_1\to S\] which we denote by $\alpha(x,f)=x \cdot f$. This action is defined for all pairs $(x,f)$ for which
$p(x)=t(f)$ and it satisfies the following axioms:
\[p(x \cdot f)=s(f), \hskip0.5cm (x \cdot f) \cdot g = x \cdot (fg), \hskip0.5cm x \cdot id_{p(x)} = x.\]
A map between $\BC$-sheaves or a {\it $\BC$-equivariant map} is a map of \'etale spaces over $C_0$ which is compatible with the
$\BC$-action.
\end{definition}

\begin{proposition} The category of $\BC$-sheaves is a topos.
\end{proposition}
\begin{definition} The topos ${{\mathcal B}\BC}$ of $\BC$-sheaves is called \bf{the classifying topos of the s-\'etale topological category}
$\BC$.
\end{definition}

\begin{examples} We provide now some examples of $\BC$-sheaves to illustrate their significance:

(i) Any small category $\BC$ can be seen as a topological category with the discrete topology. Then a $\BC$-sheaf is the same thing as a
presheaf on $\BC$ which justifies the same notation ${{\mathcal B}\BC}$ as in the Definition 5.2 for the classifying topos of a small category
$\BC$.

(ii) Any topological space $X$ may be seen as a discrete topological category ${\mathbb X}$ (the one for which all morphisms are identities).
Then a ${\mathbb X}$-sheaf is just a sheaf on $X$ and the topos ${\mathcal B}{\mathbb X}$ is the topos $Sh(X)$ of sheaves on $X$.

(iii) Let ${\mathcal G}$ be an action groupoid coming from the right action of a topological group $G$ on a topological space $X$. The groupoid
${\mathcal G}$ has $X$ as space of objects and $X \times G$ as space of morphisms where morphisms are of the form $(x,g) \colon x \cdot g \to
x$. Then a ${\mathcal G}$-sheaf $p\colon S \to X$ is a sheaf which is $G$-equivariant. Therefore ${\mathcal B}{\mathcal G}$ is the category of
$G$-equivariant sheaves.

\end{examples}

In the case of an s-\'etale topological category, i.e. a topological category with the source map $s\colon C_1 \to C_0$ being an \'etale map we
have the following definition.

\begin{definition} Let $\BC$ be an s-\'etale topological category.
A {\bf $\BC$-bundle over a space $X$} is an \'etale map (sheaf) $p\colon E\to X$ with a continuous fibrewise left action given by maps
\[\pi\colon E\to B_0, \hskip0.5cm {\rm and} \hskip0.5cm a\colon
B_1\times_{B_0} E\to E.\] Such a $\BC$-bundle is called {\bf principal} if the three conditions of non-emptiness, transitivity and freeness hold
for each $x \in X$:

(i) The stalk $E_x$ is non-empty.

(ii) For any two points $y\in E_x$ and $z\in E_x$, there are a $w \in E_x$ and arrows $\alpha\colon \pi(w) \to \pi(y)$ and $\beta\colon
\pi(w)\to \pi(z)$ such that $\alpha w=y$ and $ \beta w=z$.

(iii) For a any point $y\in E_x$ and any pair of arrows $\alpha, \beta$ in $\B$ with $s(\alpha)=\pi(y)=s(\beta)$ and $\alpha y = \beta y$ there
is a point $w\in E_x$ and an arrow $\gamma\colon \pi(w)\to \pi(y)$ in $\BC$  such that $\gamma w= y$ in $E_x$ and $\alpha \gamma = \beta \gamma$
and $\gamma z = y$ in $\B$.

A map between two principal $\BC$-bundles is a sheaf map preserving the $\BC$-action. The resulting category of principal $\BC$-bundles will
again be denoted as ${\rm Prin}(X,\BC)$.
\end{definition}

\begin{remark} \label{equivDefTop} A small category can be viewed as an s-\'etale
topological category with the discrete topology. In this case the respective definitions of principal bundles and of classifying topoi are of
course equivalent. A topological category is locally connected if the spaces of objects and arrows are locally connected. For a locally
connected s-\'etale topological category the classifying topos introduced in this section and the one defined as the topos of sheaves on the
nerve are weak homotopy equivalent.
\end{remark}
\vskip0.5cm In both cases (small and s-\'etale topological) we have the same notion of concordance of principal $\BC$-bundles as in topological
case (see \ref{conc}).

\vskip0.5cm For either a small or an s-\'etale topological category we have:

\begin{theorem}\label{equive} There is a natural equivalence of
categories
\[{\Hom(Sh(X),{\mathcal B} \BC})\simeq {\rm Prin}(X,\BC).\] On homotopy classes of topos morphisms we
have the natural bijection
\[[Sh(X), {\mathcal B} \BC]\cong {\rm Prin}_c(X,\BC).\]
For a CW complex X and any small category or any locally contractible s-\'etale category  $\BC$  there is a natural bijection
\[[X, B\BC] \cong {\rm Prin}_c(X,\BC),\] where in s-\'etale case the
fat geometric realization is taken in order to construct the classifying space.
\end{theorem}

\begin{proposition} For either a small category or a locally
connected s-\'etale category there is a natural weak homotopy equivalence
\[Sh(B\BC) \to {\mathcal B}\BC.\]
\end{proposition}

\begin{remark} The Definition \ref{Sh(Y)} of the topos $Sh(Y)$ of sheaves on the simplicial space $Y$
generalizes to the case when the opposite simplicial model category $\Delta^{\rm op}$ is replaced by an arbitrary small category $\BK$. Then,
instead of a simplicial space, we have a diagram of spaces indexed by $\BK$, i.e. a covariant functor $Y$ from $\BK$ into the category $Top$
topological spaces. With an evident modification of the Definition \ref{simpsp} we obtain the topos of sheaves on the diagram of spaces $Y$.
\end{remark}

\begin{remark}\label{Groth} From a diagram of spaces indexed by a small category $\BK$
we can construct a category $Y_{\BK}$. Object is a pair $(k,y),\, k\in \BK,\, y \in Y_k$ and arrow $(k,y)\to (l,z)$ is an arrow in $\BK$
$\alpha\colon k\to l$ such that $Y(\alpha)(y)=z$.  This is just the Grothendieck construction. The category $Y_{\BK}$ can be equipped with an
s-\'etale topology. Further, a diagram of spaces $Y$ is called locally contractible if each $Y_k$ is locally contractible. For a locally
contractible $Y$ the Grothendieck construction gives a locally contractible s-\'etale topological category $Y_{\BK}$.
\end{remark}

\begin{proposition} \label{equiv11} Let $Sh(Y)$ be the category of sheaves on a diagram of
spaces $Y$ indexed by a small category $\BK$. Then there is a natural equivalence of topoi
\[Sh(Y)\simeq {\mathcal B}(Y_{\BK}).\] Hence, for any topological space $X$
there is a natural equivalence
\[Hom(Sh(X), Sh(Y))\simeq {\rm Prin}(X, Y_{\BK}).\]
\end{proposition}
\vskip0.5cm A Principal $Y_{\BK}$-bundle canalso  be characterized as a principal $\BK$-bundle equipped with an augmentation. Let us recall that
a principal $\BK$-bundle over $X$ consists of a system of sheaves $E^k$ for each object $k$ of $\BK$ on $X$ and sheaf maps $E(\alpha) \colon E^k
\to E^l$ for each arrow $\alpha \colon k \to l$.  An augmentation on of $E$ over $Y$ is a system of maps ${\rm aug}^k \colon E^k \to Y_k$ such
that for any arrow $\alpha \colon k\to l$
\[Y(\alpha){\rm aug}^k = {\rm aug}^l E(\alpha).\]
Together with morphisms of principal bundles which respect augmentations we have the category
\[{\rm AugPrin}(X,\BK,Y)\] of principal $\BK$-bundles with an
augmentation to $Y$.

\begin{proposition} For $X$ and $Y$ as above we have a natural
equivalence of categories
\[Hom(Sh(X), Sh(Y)) \simeq {\rm Prin}(X, Y_{\BK}) \simeq {\rm AugPrin}(X,\BK,Y).\]
\end{proposition}

\begin{remark} \label{si} The case $\BK = \Delta^{\rm op}$ gives the Theorem
\ref{equiv1} as a Corollary. For this, the following equivalence
\[{\rm Prin}(X, {\Delta^{\rm op}})\simeq Lin(X)\] has to be used.
A principal $\Delta^{\rm op}$-bundle $E$ over $X$ is a simplicial sheaf such that each stalk $E_x$ is a principal $\Delta^{\rm op}$-bundle $E$
over one-point space $x$, i.e. a principal simplicial set. Finally a simplicial set is principal only if it is a nerve of a (uniquely
determined) non-empty linear order.
\end{remark}
\vskip0.5cm Next, let us consider the case $\BK = \Delta^{\rm op}\times \Delta^{\rm op}$, i.e. in this case a diagram of spaces $Y$ labeled by
$\Delta^{\rm op}\times \Delta^{\rm op}$ is just a bisimplicial space. Concerning principal $\Delta^{\rm op}\times \Delta^{\rm op}$-bundles over
$X$ we have the following result which follows from \cite{MM} (chapter VII, exercise 14).

\begin{proposition} There are natural equivalences of categories
\[{\rm Prin}(X, \Delta^{\rm op}\times \Delta^{\rm op})
\simeq {\rm Prin}(X, \Delta^{\rm op})\times {\rm Prin}(X, \Delta^{\rm op})\simeq Lin(X)\times Lin(X)\]
\end{proposition} \vskip0.5cm

Now, an augmentation is the same thing as a bisimplicial map from the product of two linear orders $NL\times NL'$ to a bisimplicial set $Y$.
Hence similarly to Definition \ref{aug} we do have
\begin{definition}\label{def} For any space $X$ and any bisimplicial space
$Y$ write $Lin^2(X,Y)$ for the product category of linear orders $(L\times L',\rm{aug})$ over $X$ equipped with a bisimplicial map ({\bf
augmentation}) ${\rm{aug}}\colon NL \times NL'\to Y$ from the product of nerves of $L$ and $L'$ to $Y$. Morphism $(L\times L',{\rm{aug})}\to
(L_1 \times L_1',\rm{aug}')$ in $Lin^2(X,Y)$ are maps of products of linear orders $L\times L'\to L_1 \times L_1'$ such that the induced map $NL
\times NL'\to NL_1\times NL_1'$ on the products of nerves respects the augmentations.
\end{definition}
\vskip0.5cm With the same definition of concordance as in \ref{conc} we do have similarly to Theorem \ref{equiv1}:
\begin{theorem}\label{equiv10}Let $Y$ be a bisimplicial space. For any space
$X$ there is a natural equivalence of categories
\[{\rm Hom}(Sh(X), Sh(Y))\simeq Lin^2(X,Y).\] On homotopy classes of topos morphisms we
have the natural bijection
\[[Sh(X), Sh(Y)]\cong Lin^2_c(X,Y).\]
\end{theorem} \vskip0.5cm Similarly to Theorem \ref{weak-equiv1} we
have the following Theorem, where the geometric realization $|Y|$ of a a bisimplicial
space $Y$ can be taken as the geometric realization of its diagonal. Equivalently $Y$ can be defined as the "horizontal" geometric realization
followed by the the "vertical" one or vice versa.

\begin{theorem} \label{weak-equiv-bi} For any bisimplicial space $Y$ the
topoi $Sh(Y)$ and $Sh(|Y|)$ have the same weak homotopy type.
\end{theorem}
\vskip0.5cm We recall that, in accordance with Remark \ref{Groth}, a bisimplicial space $Y$ is locally contractible if all spaces $Y_{n,m}$ are
locally contractible. Again, geometric realization of a locally contractible bisimplicial space is locally contractible. Hence, we have
similarly to \ref{classifi} the following corollary:
\begin{corollary}\label{classifi13} Let $X$ be a CW-complex and
$Y$ be a locally contractible bisimplicial space. There is a natural bijection between homotopy classes of maps $[X,|Y|]$ and concordance
classes $Lin^2_c(X,Y)$.
\end{corollary}

\section{Classifying spaces and classifying topoi of topological bicategoies II}

It is beyond the scope of this paper to give a full account of the constructions of Lack and Paoli, Tamsamani and Simpson. Concerning the latter
two  the interested reader may find useful the nice survey of definitions of $n$-categories by T. Leinster \cite{Lei}. Let $\mathbf{\rm Set}$
and $\mathbf{\rm Cat_1}$ denote the categories of (small) sets and (small) categories, respectively, and let $\mathbf{\rm Cat}$ denote the
2-category of (small) categories.

\vskip0.5cm
\subsection{Lack-Paoli nerve as a singular functor} The nerve
construction of Lack and Paoli \cite{LP} is obtained as the singular functor when ${\mathcal V}=\mathbf{\rm Cat}$. In order to define the nerve
$N\B$ of a (small) bicategory $\B$ they introduced a (strict) 2-category $\mathbf{\rm NHom}$ with bicategories as objects, whose 1-cells are
normal homomorphisms (normal lax functors with invertible comparison maps). We will not give the general definition of 2-cells (icons) here. We
describe them below explicitly in a special case.

Every category (in particular $\mathbf{\rm Cat_1}$ and the category defined by the ordinal [n]) can be seen as a locally discrete bicategory with only identity 2-cells. The normal
homomorphism between locally discrete bicategories is just a functor between the corresponding categories, and there are no nontrivial icons
between such. In this way, we obtain a fully faithful inclusion 2-functor
\[\label{LPemb} J\colon \Delta \to \mathbf{\rm NHom}\]
and the category $\Delta$ can be seen as a full sub-2-category of $\mathbf{\rm NHom}$. The singular 2-functor
\[\label{LP2nerve}
N_{LP} \colon\mathbf{\rm NHom} \to [\Delta^{\rm op},\mathbf{\rm Cat}]\] of the inclusion $J$ is Lack and Paoli 2-nerve. The 2-functor $N_{LP}$
is fully faithful.

\begin{definition} A {\bf normal homomorphism} $(B,f,\beta)\colon [n]\to \B$ from an ordinal $[n]$ to the bicategory $\B$ is a lax normal functor
for which each 2-cell $\beta_{ijk}$ in the Definition 4.1 is invertible.
\end{definition}

\begin{definition} An {\bf icon} between normal homomorphisms $F,G \colon [n] \to \B$ of bicategories is a lax natural transformation
$\phi \colon F \Rightarrow G$, in which the component $\phi_i \colon B_i \to C_i$ is an identity, for each $i \in [n]$. More explicitly, an icon
$\phi \colon (B, f, \beta) \Rightarrow (C, g, \gamma)$ consists of the following data:
\begin{itemize}
\item [(i)] for any $i \in [n]$ an identity $B_i=C_i$

\item [(ii)] for each $i,j \in [n]$ such that $i \leq j$, a 2-cell $\phi_{ij}\colon f_{ij} \Rightarrow
g_{ij}$
\[\xymatrix@!=4pc{B_j \rtwocell<8>^{f_{ij}}_{g_{ij}}{\,\,\,\,\,\,\phi_{ij}} & B_i}\]
such that  for all $i, j, k\in [n]$ with $i \leq j \leq k$ we have an equality of pasting diagrams
\[\xymatrix@!=3pc{& B_j \ar@/^1.5pc/[dr]^{f_{ij}} &\\
B_k \rrtwocell<9>^{f_{ik}}_{g_{ik}}{\,\,\,\,\,\,\phi_{ik}} \rrtwocell<\omit>{<-8>\,\,\,\,\,\,\beta_{ijk}} \ar@/^1.5pc/[ur]^{f_{jk}} && B_i}
\hspace{2cm} \xymatrix@!=3pc{&\\=}
\xymatrix@!=3pc{& B_j \drtwocell<6>^{f_{ij}}_{g_{ij}}{\,\,\,\,\,\,\phi_{ij}} &\\
B_k \ar@/_2pc/[rr]_{g_{ik}} \rrtwocell<\omit>{\,\,\,\,\,\,\gamma_{ijk}} \urtwocell<6>^{f_{jk}}_{\,\,\,\,\,\,g_{jk}}{\,\,\,\,\,\,\phi_{jk}} &&
B_i}\] which means that the following identity of 2-cells holds in $\B$:
\[\label{iconcoh}
\phi_{ik}\beta_{ijk}=\gamma_{ijk}(\phi_{ij} \circ \phi_{jk}).
\]
\end{itemize}
\end{definition}
\subsection{Characterization of Lack-Paoli 2-nerves of bicategories}\label{condi}
In their paper \cite{LP}, Lack and Paoli also described necessary and sufficient conditions for a simplicial object $X\colon \Delta^{\rm op} \to
{\mathbf{\rm Cat}}$ to be a 2-nerve of a bicategory. In order to provide such characterization, they used {\it discrete isofibrations} which are
functors $P\colon E \to B$ such that for each object $e$ in the category $E$ and each isomorphism $\beta\colon b\to P(e)$ in $B$ there exists a
unique isomorphism $\varepsilon \colon e' \to e$ in $E$ with $P(\varepsilon) = \beta$. Let further ${\rm {c}}_n \colon X_n \to {\rm
Cosk}_{n-1}(X)_n$ denotes the $n$-component of the simplicial map $c \colon X \to {\rm Cosk}_{n-1}(X)$ from a simplicial object $X$ to its
{n-1}-coskeleton ${\rm Cosk}_{n-1}(X)$, which is the unit of an adjunction between (n-1)-truncation $tr_n$ and (n-1)-coskeleton ${\rm
Cosk}_{n-1}$.

\begin{theorem} The necessary and sufficient conditions for a 2-functor $X\colon
\Delta^{\rm op} \to {\mathbf{\rm Cat}}$ to be a 2-nerve of a bicategory are:
\begin{itemize}
\item [(i)] $X$ is 3-coskeletal,

\item [(ii)] $X_0$ is discrete,

\item [(iii)] the Segal functors $S_n\colon X_n \to X_1 \times_{X_0}\ldots
\times_{X_0}X_1$ are equivalences of categories,

\item [(iv)] ${\rm {c}}_2$ and ${\rm {c}}_3$ are discrete isofibrations.
\end{itemize}
\end{theorem}
\subsection{Lack-Paoli 2-nerve as a bisimplicial set (space)} \label{2-nerve bisimpl}
If we apply the Grothendieck nerve functor at each level of the 2-nerve of Lack and Paoli (\ref{LP2nerve}), we obtain a functor
\[\label{BSnerve} B_{LP} \colon\mathbf{\rm NHom} \to [\Delta^{\rm
op},\mathbf{\rm SSet}]
\]
where the right-hand side is the category of bisimplicial sets. If we define the 2-nerve in such bisimplicial terms, the definition makes sense also
for a topological bicategory $\B$, in which case the 2-nerve will  naturally be a bisimplicial space. Although the above conditions \ref{condi}
can be translated into the bisimplicial language, we will not do it here. From now on we will understand the 2-nerve of Lack \& Paoli as a
bisimplicial set (bisimplicial space in case of a topological bicategory).

\subsection{Tamsamani and Simpson}
Let ${\mathbf{\rm Tam}}$ denote the full sub-2-category of $[\Delta^{\rm op},\mathbf{\rm Cat}]$ consisting of those $X$, for which $X_0$ is
discrete and the Segal maps $S_n$ are equivalences. Let further ${\mathbf{\rm Simpson}}$ denote the smaller full sub-2-category of those $X$ for
which the Segal maps $S_n$ are fully faithful and surjective on objects. Also, in these cases we can interpret these "2-nerves" as bisimplicial
sets (see \cite{T}, \cite{S}, \cite{Lei}, where the corresponding definitions  can be found). We will speak of Tamsamani 2-nerve (or 2-category)
and Simpson 2-nerve (or 2-category).

\begin{remark} \label{6.8} The Lack-Paoli  2-nerve is also a Simpson 2-nerve
and thus also a Tamsamani 2-nerve. To each Tamsamani 2-nerve $X$ there is a bicategory $GX$  (and vice versa) constructed in \cite{T}. We refer
the reader for more details on Tamsamani 2-nerves (including a proper notion of equivalence) to this paper.
\end{remark} Here we only mention the following results of Lack \&
Paoli: \vskip0.5cm The (Lack-Paoli) 2-nerve 2-functor $N_{LP}\colon\mathbf{\rm NHom} \to {\mathbf{\rm Tam}}$, seen as landing in the 2-category
${\mathbf{\rm Tam}}$, has a left 2-adjoint $G$. Since $N_{LP}$ is fully faithful, the counit $GN_{LP} \to 1$ is invertible. Each component
$u\colon X \to N_{LP}G$ of the unit is a pointwise equivalence (i.e. each component $u_n$ is an equivalence) and $u_0$ and $u_1$ are identities.
\vskip0.5cm Let $\mathbf{\rm Ps}(\Delta^{\rm op}, \mathbf{\rm Cat})$ denotes the 2-category of 2-functors, pseudonatural transformations and
modifications and let ${\mathbf{\rm Tam_{ps}}}$ be its full sub-2-category consisting of Tamsamani 2-categories. Then the 2-nerve 2-functor
$N_{LP} \colon\mathbf{\rm NHom} \to {\mathbf{\rm Tam_{ps}}}$ is a biequivalence of 2-categories.

\begin{definition}\label{cltop}Let $\B$ be a topological bicategory. The {\bf classifying topos} ${\mathcal B}_{LP}\B$ of the topological
bicategory $\B$ is defined as the topos of sheaves $Sh(B_{LP}\B)$ on the bisimplicial space $B_{LP}\B$ (Lack-Paoli bisimplicial nerve).
\end{definition}
\begin{definition}\label{clsp}The {\bf classifying space} $B_{LP}\B$ of a topological bicategory
$\B$ is the geometric realization $|B_{LP}\B|$ of its bisimplicial nerve $B_{LP}\B$.
\end{definition}
\vskip0.5cm With these definitions we have the following corollary of Theorem \ref{weak-equiv-bi}.
\begin{corollary}\label{cor6.11}For any topological bicategory $\B$ the topos
of sheaves $Sh(B_{LP}\B)$ on the classifying space $B_{LP}\B$ has the same weak homotopy type as the classifying topos ${\mathcal B\mathcal
B}\B$.
\end{corollary}
\begin{definition}\label{Lin2}For a topological bicategory $\B$
write $Lin^2(X,\B)$ for the product category of linear orders over $X$ equipped with an augmentation ${\rm{aug}}\colon NL \times NL'\to NN\B$.
An object $E$ of this category will be called a {\bf Lack-Paoli principal $\B$-bundle}. We call two Lack-Paoli principal $\B$-bundles $E_0$ and
$E_1$ on $X$ {\bf concordant} if there exists a Lack-Paoli principal $\B$-bundle on $X\times [0,1]$ such that we have the equivalences
$E_0\simeq i^*_0(E)$ and $E_0\simeq i^*_1(E)$ under the obvious inclusions $i_0, i_1\colon X \hookrightarrow X\times [0,1]$.
\end{definition} \vskip0.5cm Similarly to Theorems
\ref{equiv2} and \ref{equivDus} we have from Theorem \ref{equiv10} the following "classifying" property of the classifying topos ${\mathcal
B}_{LP}\B.$
\begin{theorem}\label{th6.13}For a topological bicategory $\B$ and a
topological space $X$ there is a natural equivalence of categories
\[{\Hom (Sh(X),{\mathcal B}_{LP} \B})\simeq Lin^2(X,\B).\] On homotopy classes of topos morphisms we
have the natural bijection
\[[Sh(X), {\mathcal B}_{LP} \B]\cong Lin^2_c(X,\B).\]
\end{theorem}
\vskip0.5cm The "classification" property of the classifying space $B_{LP}\B$ now follows as a corollary from Corollary \ref{classifi13}.
\begin{corollary}\label{cor6.14}For a locally contractible bicategory
$\B$ and a CW-complex $X$ there is a natural bijection
\[[X, B_{LP}\B] \cong Lin^2_c(X,\B).\]
\end{corollary}
\begin{remark} ({\bf Tamsamani and Simpson principal $\B$-bundles})
In the above Definitions  \ref{cltop}, \ref{clsp} and \ref{Lin2} we could have used instead of Lack-Paoli 2-nerve the Tamsamani or Simpson
2-nerve (in the case these are bisimplicial spaces). Obviously, for such {\it Tamsamani} and {\it Simpson} principal $\B$-bundles Corollaries
\ref{cor6.11}, \ref{cor6.14} and Theorem \ref{th6.13} are still valid.
\end{remark}


\end{document}